# Geometry in the large on Hadamard manifolds

Stepanov Sergey[1,2], Tsyganok Irina[2]

[1] Department of Mathematics, Russian Institute for Scientific and Technical Information of the Russian Academy of Sciences, 20, Usievicha str., A-190, Moscow, Russia
s.e.stepanov@mail.ru

[2] Department of Mathematics, Financial University, 49, Leningradsky pr., Moscow 14993, Russia
e-mail: i.i.tsyganok@mail.ru

**Abstract**. In this paper, we prove several Liouville-type theorems on the non-existence of Killing-Yano tensors, Killing tensors, and harmonic symmetric tensors on Hadamard manifolds and, in particular, on Riemannian symmetric spaces of non-compact type. These theorems supplement the well-known vanishing theorems for the above tensors, obtained using the Bochner technique for compact Riemannian manifolds. In turn, the proofs of our theorems will use well-known Liouville-type theorems on the non-existence of subharmonic and harmonic functions on complete Riemannian manifolds, which we have partially modified for the case of Hadamard manifolds and, in particular, Riemannian symmetric spaces of noncompact type.

**Mathematics Subject Classification (2010):** 53C20

## 1. Subharmonic, harmonic and convex $L^q$-functions on Hadamard manifolds

The most rigid classical concept of curvature, and therefore containing most of the information about a Riemannian manifold, is sectional curvature. One aspect that makes the concept of sectional curvature particularly important is its relationship to the topology of a Riemannian manifold. For example, as a starting point in the study of the geometry of Riemannian manifolds of non-positive sectional curvature, any geometer will first recall the well-known *Cartan-Hadamard theorem*: Let $(M, g)$ be an $n$-dimensional $(n \geq 2)$ simply connected and complete Riemannian manifold of non-positive sectional curvature, then $(M, g)$ is diffeomorphic to the Euclidean space $\mathbf{R}^n$. Therefore, a simply connected complete Riemannian manifold of non-positive curvature was called a *Hadamard manifold* after the Cartan-Hadamard theorem (see, for example, [1, p. 241]; [2, pp. 391-381]; [5]). Well-known basic examples of such manifolds are the Euclidean space $\mathbf{R}^n$, with zero sectional curvature, and the hyperbolic space $\mathbf{H}^n$, with constant negative sectional curvature. On the other hand, a flat torus $\mathbf{T}^n$ is a complete manifold of zero sectional curvature. At the same time, $\mathbf{T}^n$ is connected, but it is not a simply connected manifold and therefore is not an example of Hadamard manifolds. Therefore, the requirement of simply connectedness is essential here and, for example, strong enough to distinguish between Euclidean space $\mathbf{R}^n$ and a flat torus $\mathbf{T}^n$.

Side by side, we can formulate two obvious corollaries follow from the Cartan-Hadamard theorem. First, one can conclude from the theorem that no compact simply connected manifold admits a metric of non-positive sectional curvature (see [1, p. 162]). Second, it also follows from the theorem that a Hadamard manifold has an infinite volume (see [4, p. 4732]).

The next part of the section will be devoted to the relationships between the geometry of a Riemannian manifold and the global behavior of its subharmonic and harmonic functions under assumptions on sectional curvature (see, for, example, [1, p. 281-284]; [2]; [8-10] and etc.). Most of these results are called *Liouville-type theorems* or, in other words, *vanishing theorems* and belong to the *Bochner technique* (see, for example, [1. p. 333-364]; [3]; [7] and etc.). We recall here that a scalar function $f \in C^2(M)$ is *subharmonic* (see [1, p. 281]) if it is satisfies the differential inequality $\Delta f \geq 0$ for the *Beltrami Laplacian* $\Delta = div \circ grad$ and, in particular, $f \in C^2(M)$ is *harmonic* (see [1, p. 283]) if it is a smooth solution of the Laplace equation $\Delta f = 0$. It is well-known that in mathematics, subharmonic and harmonic functions are important classes of functions used extensively in partial differential equations, complex analysis and potential theory. These functions are intensively studied in the geometry of Riemannian manifolds. First, the well-known Hopf lemma (see [11, p. 338]) shows that a compact Riemannian manifold has no subharmonic and harmonic functions except for constant functions. On the other hand, Huber (see [12]) proved that a complete two-dimensional Riemannian manifold with non-negative curvature does not admit a non constant negative subharmonic function. Next, Karp (see [13]) found that a complete noncompact Riemannian manifold admits no non constant negative subharmonic function if it has moderate volume growth. For comparison, we recall that, the classical Liouville theorem says that a subharmonic function defined over $\mathbf{R}^2$ (or a harmonic function defined over $\mathbf{R}^n$) and bounded above is constant. Several further results about the properties of subharmonic and harmonic functions on complete Riemannian manifolds were obtained by many authors, such as Li and Schoen, Greene and Wu, Yau and others (see, for example, [3]; [8]).

The well-known question is still very important: which subspaces of subharmonic (and harmonic) functions on complete Riemannian manifolds contain only constant functions and, in particular, identically equal to zero. For this case the follows Yau theorem is well-known (see [9]): If $(M, g)$ is a complete noncompact Riemannian manifold (without boundary) and $f : M \to [0, +\infty)$ is a subharmonic $L^q$-function for $1 < q < +\infty$, then $f$ is a constant.

Throughout this article, we will consider scalar functions from $L^q(M)$ where $0 < q < +\infty$ and $(M, g)$ is a complete noncompact Riemannian manifold without boundary. By definition, each such scalar function $f$ defined on $(M, g)$ and belonging to $L^q(M)$ for $0 < q < \infty$ must satisfy the integral inequality $\int_M |f|^q dvol_g < \infty$ where the integral will always be understood in terms of the volume form $dvol_g$ of the metric $g$. In addition, here the situation divides into two cases. The first case is when $(M, g)$ has finite volume, then all constant functions are in $L^q(M)$ for any $0<q<\infty$. The second case is when $(M,g)$ has infinite volume, then all of the constant functions, but zero, are in $L^q(M)$ for any $0<q<\infty$. Namely, if the function $f \in L^q(M)$ for any $q \in (0,\infty)$ is a constant $C$, then the inequality $\int_M |f|^q dvol_g < \infty$ becomes $|C|^q \cdot \int_M dvol_g < +\infty$. If, in addition, $(M,g)$ has an infinite volume, then from the last inequality we obtain $C=0$.

In turn, using the above information about Hadamard manifolds, we can conclude that the following lemma holds.

**Lemma 1.1.** *The Hadamard manifold $(M,g)$ does not admit a non-zero non-negative subharmonic $L^q$-function for each $q \in (0,+\infty)$.*

*Proof.* Let $(M,g)$ be a Hadamard manifold and $f \in C^2(M)$ be a non-negative subharmonic $L^q$-function for $q \in (0,+\infty)$. At the same time, it is known from [10] that on a complete simply connected Riemannian manifold $(M,g)$ of non-positive sectional curvature every non-negative subharmonic function $L^q$-function for any $q \in (0,\infty)$ is a constant $C$. We recall that Hadamard manifolds have infinite volume, and hence $C=0$. This finishes the proof of our lemma.

The well-known fundamental fact is that if $f \in C^2(M)$ is a harmonic function, then $|f|^q$ is a nonnegative subharmonic function for each $q \geq 1$ (see also [3, p. 373]). Therefore, if every nonnegative subharmonic $L^q$-function on $(M,g)$ is constant, then every harmonic $L^q$-function on $(M,g)$ is constant too (see [10]). Therefore, we can formulate a corollary from Lemma 1.1.

**Corollary 1.2**. *The Hadamard manifold $(M,g)$ does not admit a non-zero harmonic $L^q$-function for $q \in (0,+\infty)$.*

Next, we formulate a Liouville-type proposition on a smooth $L^q$-function on a Hadamard manifold.

**Lemma 1.3.** *The Hadamard manifold $(M,g)$ does not admit a non-zero smooth nonnegative $L^q$-function $f \in C^2(M)$ for each $q \in (0, +\infty)$ and $q \neq 1$ such that $f \cdot \Delta f \geq 0$.*

*Proof.* Let $(M,g)$ be a Hadamard manifold. It is known from [9] that if the inequality $f \cdot \Delta f \geq 0$ holds for a nonnegative function $f \in C^2(M)$ defined on a complete Riemannian manifold $(M,g)$, then either $\int_M f^q dvol_g = +\infty$ for all $q \neq 1$ or $f = constant$. In particular, here $q$ may even be less than one (see also [9]). In this case, the inequality $\int_M f^q dvol_g < \infty$ becomes $C^q \cdot \int_M dvol_g < \infty$. At the same time, we known that a Hadamard manifold $(M,g)$ have infinite volume, and hence from the last inequality we obtain $C = 0$. This completes our proof.

The third part of the section, we consider the concept of convex functions. Recall that $f \in C^2(M)$ is *convex function* if, at each point, its Hessian $Hess_g f := \nabla^2 f$ is positive semi-definite (see [14]). In this case, we have $\Delta f \geq 0$ and hence $f$ is a subharmonic function. In differential geometry the existence of convex functions on a $(M,g)$ is a long standing problem. The first solution of this problem can be found in the article [14] of Bishop and O'Neill. They proved that if $(M,g)$ is complete and has finite volume, then it does not possess a non-constant smooth convex function. We supplement their result on the basis of Corollary 1.2 by formulating the following proposition.

**Corollary 1.4**. *A Hadamard manifold does not admit non-zero non-negative convex $L^q$-functions for some $q \in (0, +\infty)$.*

Moreover, the article [14] considers some general methods for constructing convex functions on a Hadamard manifold $(M,g)$. Consider an example constructed using the Riemannian distance function $d(x,\cdot)$ on $(M,g)$. We recall that for an isometry $\phi:(M,g)\to(M,g)$ its *displacement function* $f:(M,g)\to(0,+\infty)$ is given by the equality $f(x)=[d(x,\phi(x))]^2$. In the case of a Hadamard manifold $(M,g)$, the displacement function $f$ must be non-negative convex and smooth (see [14]; [15]). Therefore, the corollary is true.

**Corollary 1.5**. *Let $(M,g)$ be a Hadamard manifold and $f$ be the displacement function $f(x)=[d(x,\phi(x))]^2$ of an isometry $\phi:(M,g)\to(M,g)$. If $f$ is a $L^q$-function for some $q\in(0,+\infty)$, then $f\equiv 0$ and hence $\phi$ is the identity.*

The fourth part of the section will be devoted to finite-dimensional Riemannian symmetric spaces. The condition of the parallelism of the curvature tensor $\nabla R=0$ defines the class of *Riemannian locally symmetric spaces* that can be equivalently defined as those Riemannian manifolds which are locally reflectionally geodesic symmetric around any point $x\in M$ (see [11, pp. 243-244]). A Riemannian locally symmetric space $(M,g)$ is called a *Riemannian globally symmetric space* if it's locally geodesic symmetries are defined on the whole of it. In this case, a Riemannian symmetric space $(M,g)$ is complete (see [11, p. 244]; [16, p. 240]). Riemannian globally symmetric spaces can be classified using their isometry groups. The classification distinguishes three basic types of such spaces: spaces of so-called compact type, spaces of so-called noncompact type and spaces of Euclidean type (see, e.g., [11, p. 42]; [16, p. 245]). In particular, a Riemannian globally symmetric space $(M,g)$ of non-compact type is simply connected (and therefore diffeomorphic to $\mathbf{R}^n$) and has non-positive sectional curvature (see [16, pp. 245; 246]). Therefore, a Riemannian globally symmetric space of non-compact type has an infinite volume. Using the above, we can state that Riemannian globally symmetric space of non-compact type is a prominent example of a Hadamard manifold. Therefore, a new corollary of Lemma 1.1 holds.

**Corollary 1.6**. *A Riemannian globally symmetric space of noncompact type does not admit a nonzero non-negative subharmonic $L^q$-function for $q \in (0, \infty)$.*

**Remark**. The geometry of a Riemannian symmetric space with nonpositive sectional curvature is described in detail in [24]. In addition to the above, we note that Hopf's lemma (see [11, p. 338]) shows that a Riemannian globally symmetric space of compact type has no subharmonic and harmonic functions, except for constant functions.

## 2. The Sampson Laplacian on Hadamard manifolds

Let $C^\infty(S^p M)$ be the space of $C^\infty$-sections of the bundle $S^p M = S^p(T^*M)$ of covariant symmetric $p$-tensors on a connected Riemannian manifold $(M, g)$. Then it is obvious that the following equality is true

$$\dim S^p(T_x^* M) = \binom{n+p-1}{p}.$$

for the vector space $S^p(T_x^* M)$ of covariant symmetric $p$-tensors on $T_x M$ at an arbitrary point $x \in M$.

Now, we define the differential operator $\delta^* : C^\infty S^p M \to C^\infty S^{p+1} M$ of degree one by the formula $\delta^* \varphi = (p+1) S^{p+1}(\nabla \varphi)$ for an arbitrary $\varphi \in C^\infty(S^p M)$ and the standard point-wise symmetry operator $S^{p+1} : T^*M \otimes S^p(T^*M) \to S^{p+1}(T^*M)$. There exists its formal adjoint operator $\delta : C^\infty S^{p+1} M \to C^\infty S^p M$ for $\delta^*$ (see [17, pp. 34-35]). Using the above, Sampson determined in [18, p. 147] the operator

$$\Delta_S = \delta \delta^* - \delta^* \delta : C^\infty S^p M \to C^\infty S^p M \qquad (2.1)$$

for an arbitrary Riemannian manifold $(M, g)$. Moreover, he showed that the operator $\Delta_S$ admits the Weitzenböck decomposition (see also [19])

$$\Delta_S = \overline{\Delta} - \mathfrak{R} \tag{2.2}$$

where $\overline{\Delta} = \delta \nabla$ is the *Bochner Laplacian* and $\mathfrak{R}$ is the *Weitzenböck curvature operator* of the ordinary Lichnerowicz Laplacian $\Delta_L = \overline{\Delta} + \mathfrak{R}$ (see [17, p. 54]; [70, p. 315]) which was restricted to covariant symmetric $p$-tensors. The Weitzenböck curvature operator $\mathfrak{R}$ of the Lichnerowicz Laplacian $\Delta_L$ can be algebraically (even linearly) expressed through the curvature $R$ and Ricci $Ric$ tensors of $(M, g)$. Moreover, it satisfies the identities (see [70, p. 315])

$$g(\mathfrak{R}(T), T') = g(T, \mathfrak{R}(T')) \tag{2.3}$$

and

$$\text{trace}_g \mathfrak{R}(T) = \mathfrak{R}(\text{trace}_g T) \tag{2.4}$$

for any $T, T' \in \otimes^p T^*M$. In particular, from (2.3) any one conclude that $\mathfrak{R}_x : \otimes^p T_x^* M \to \otimes^p T_x^* M$ is a symmetric endomorphism at any point $x \in M$. For example, if $p = 2$ then we have the equality (see [17, pp. 64; 356]; [19])

$$\mathfrak{R}(\varphi) = -2 g^{km} g^{lt} R_{mijt} \varphi_{kl} + g^{kl} R_{ki} \varphi_{lj} + g^{kl} R_{kj} \varphi_{li} \tag{2.5}$$

where $\varphi_{ij}$, $R_{ijkl}$ and $R_{ij}$ denote the local components of $\varphi \in C^\infty S^2 M$, the Riemannian curvature tensor $Rm$ and the Ricci tensor $Ric$, respectively. In addition, $g^{kl}$ are the local contravariant components of the metric tensor $g$. These local components are defined by the following identities $\varphi_{ij} = \varphi(e_i, e_j)$, $R_{ijkl} = g_{im} R^m_{jkl}$ and $R_{kl} = R^i_{kil}$ where $R(e_j, e_l) e_k = R^i{}_{kjl} e_i$ and $g_{im} = g(e_i, e_m)$ for any frame $e_1, \ldots, e_n$ of $T_x M$ at an arbitrary point $x \in M$ and for any $i, j, k, \ldots = 1, 2, \ldots, n$.

**Remark.** The Sampson operator $\Delta_S$ is a Laplacian. Therefore, the kernel $\text{Ker} \Delta_S$ of $\Delta_S$ is a finite-dimensional vector space on a compact manifold $(M, g)$. More information about the properties and applications of $\Delta_S$ can be found in papers from the following list: [19]; [21-23].

In accordance with the general theory we define the two vector spaces (see [8, p. 104]). First, we define by the condition

$$\mathrm{H}(S^p M) = \{\varphi \in C^\infty(S^p M): \Delta_S \varphi = 0\}$$

the vector space of $\Delta_S$-*harmonic symmetric p-tensors* $\varphi \in C^\infty(S^p M)$. Second, we define by the condition

$$L^q \mathrm{H}(S^p M) = \{\varphi \in \mathrm{H}(S^p M): \|\varphi\| \in L^q(M) \text{ and } q \in (0, \infty)\}$$

the vector space of $\Delta_S$-*harmonic symmetric $L^q$-tensors*. The norm of $\varphi \in C^\infty(S^p M)$ with respect to the Riemannian metric $g$ is denoted by the symbol $\|\cdot\|$. Using these notations, we can conclude from (2.2) and (2.4) that if $\varphi \in \mathrm{H}(S^p M)$, then $\text{trace}_g \varphi \in \mathrm{H}(S^{p-2} M)$. Moreover, we can prove the following theorem.

**Theorem 2.1**. *Let $(M, g)$ be a Hadamard manifold then the vector space $L^q \mathrm{H}(S^2 M)$ for some $q \in (0, +\infty)$ and $q \neq 1$ is trivial.*

*Proof.* We define the non-negative scalar function $f = \|\varphi\|$ for $\varphi \in \mathrm{H}(S^p M)$. In this case, from (2.2) we deduce the well-known *Bochner-Weitzenböck formula*

$$\frac{1}{2} \Delta f^2 = \|\nabla \varphi\|^2 - g(\mathfrak{R}(\varphi), \varphi) \tag{2.6}$$

for an arbitrary $\varphi \in \mathrm{H}(S^p M)$. On the other hand, we have

$$\frac{1}{2} \Delta f^2 = f \cdot \Delta f + \|\nabla f\|^2. \tag{2.7}$$

Then from (2.7) and the Bochner-Weitzenböck formula (2.6) we obtain

$$f \cdot \Delta f = \|\nabla \varphi\|^2 - g(\Re(\varphi),\varphi) - \|\nabla f\|^2 \tag{2.8}$$

where $\|\nabla \varphi\|^2 \geq \|\nabla f\|^2$ due to the Kato inequality $\|\nabla \varphi\| \geq \|\nabla\|\varphi\|\|$ (see [7, p. 380]). In this case, from (2.8) one can obtain the inequality

$$f \cdot \Delta f \geq -g(\Re(\varphi),\varphi) \tag{2.9}$$

for an arbitrary $\varphi \in H(S^p M)$. In particular, if $p = 2$ then from (2.5) in accordance with [17, p. 436] we get the formula

$$g(\Re(\varphi),\varphi) = \sum_{i<j} \sec(e_i,e_j)(\mu_i - \mu_j)^2 \tag{2.10}$$

where $e_1, e_2, \ldots, e_n$ is the orthonormal basis of $T_x M$ at an arbitrary $x \in M$ such that $\varphi(e_i, e_j) = \delta_{ij}\mu_i$ and for the sectional curvature $\sec(e_i, e_j)$ in the direction of the two-plane $\sigma(x) = span\{e_i, e_j\}$. From (2.10) we conclude that if the section curvature of $(M,g)$ is non-positive, then $g(\Re(\varphi),\varphi) \leq 0$ for any $\varphi \in S^2 M$. Then we can conclude from (2.9) that $f \cdot \Delta f \geq 0$ where $f = \|\varphi\|$ for $\varphi \in H(S^2 M)$. Next we can refer Lemma 1.3. This completes our proof.

From the above we known that $g(\Re(\varphi),\varphi) \leq 0$ for any $\varphi \in S^2 M$ if the section curvature of $(M,g)$ is non-positive. In this case the inequality $\frac{1}{2}\Delta f^2 \geq 0$ follows from (2.6). Hence $f^2$ is a subharmonic function. Then based on Lemma 1.1, we can formulate a theorem.

**Theorem 2.2**. *An $n$-dimensional Hadamard manifold does not admit a non-zero $\Delta_S$-harmonic symmetric tensor $\varphi \in C^\infty(S^2 M)$ if the square of its norm is an $L^q$-function at least one $q \in (0, +\infty)$.*

Let $S_0^p(T_x^*M) \subset S^p(T_x^*M)$ be a space of covariant symmetric $p$-tensors which are totally traceless, that is, traceless on any pair of indices at an arbitrary point $x \in M$. Then

$$\dim S_0^p(T_xM) = \binom{n+p-1}{p} - \binom{n+p-3}{n-1}.$$

The Sampson Laplacian $\Delta_S$ maps $S_0^p M$ to itself for the bundle $S_0^p M$ of traceless symmetric $p$-tensors on $(M, g)$. This property is a corollary of the identities (2.2) and (2.4). Therefore, we can conclude that $\Delta_S : C^\infty S_0^p M \to C^\infty S_0^p M$. Then, in particular, the following theorem holds.

**Theorem 2.3.** *The Hadamard manifold $(M, g)$ does not admit a non-zero symmetric $\Delta_S$-harmonic $p$-tensor $\varphi \in C^\infty(S_0^p M)$ such that the square of its norm belongs to $L^q(M)$ at least for one $q \in (0, \infty)$.*

*Proof.* Consider a non-zero symmetric $\Delta_S$-harmonic $p$-tensor $\varphi \in C^\infty(S_0^p M)$ on a Hadamard manifold $(M, g)$. In this case, the inequality (2.6) holds. In [32] and [34] was proved that the inequality $g(\Re(\varphi), \varphi) \leq 0$ holds for every $\varphi \in S_0^p(T_xM)$ at any point $x$ of a manifold $(M, g)$ with non-positive sectional curvature. For this case, we can deduce from (2.6) that $\Delta \|\varphi\|^2 \geq 0$ for a symmetric tensor $\varphi \in C^\infty(S_0^p M)$ and therefore, $\|\varphi\|^2$ is a subharmonic function. Moreover, if $\|\varphi\|^2 \in L^q(M)$ at least for one $q > 0$, then $\varphi \equiv 0$ by Theorem 2.1.

## 3. Hadamard $A$-spaces

A Hadamard manifold $(M, g)$ is diffeomorphic to Euclidean space $\mathbf{R}^n$ but, in general, it may not be isometric to $\mathbf{R}^n$. Therefore, it is quite natural to ask under which condition Cartan-Hadamard manifold is isometric to the Euclidean space $\mathbf{R}^n$. We will give one of the answers to this question in this section.

We recall that Gray introduced the notion of $A$-space in the paper [25]. An $A$-space it is a Riemannian manifold $(M, g)$ whose Ricci tensor $Ric$ satisfies the following condition: $(\nabla_X Ric)(X,X) = 0$ for all $X \in TM$. Recently, there have been many papers on the geometry of these manifolds (see [26 - 29] and others). In particular, from the classical monograph [17, p. 451] it is known that if $(M, g)$ is a compact (without boundary) $A$-space with non-positive sectional curvature, then $\nabla Ric = 0$. If in addition there exists a point in $M$ where the sectional curvature of every two-plane is strictly negative, then $(M, g)$ is Einstein, i.e. its Ricci tensor satisfies $Ric = \rho g$ for some constant $\rho$ (see [17, p. 451]). Let us complete this theorem using our Theorem 2.1. Namely, it is known from the Cartan-Hadamard theorem that an $n$-dimensional Hadamard manifold $(M, g)$ is diffeomorphic to Euclidean space $\mathbf{R}^n$. We now prove the following proposition on the condition under which an $n$-dimensional Hadamard manifold $(M, g)$ is isometric to the Euclidean space $\mathbf{R}^n$.

**Theorem 3.1.** *The Ricci tensor $Ric$ of an $n$-dimensional $A$-space $(M, g)$ belongs to the vector space $\mathrm{H}(S^2 M)$. However, if $(M, g)$ is a Hadamard $A$-space and $\|Ric\| \in L^q \mathrm{H}(S^2 M)$ at least for one $q \in (0, +\infty)$ and $q \neq 1$, then $(M, g)$ is isometric to Euclidean space $\mathbf{R}^n$.*

*Proof.* Let now $(M, g)$ be an $n$-dimensional $A$-space then its Ricci tensor $Ric$ satisfies the equations $\delta^* Ric = 0$ and has a constant trace, i.e., the scalar curvature $s = \mathrm{trace}_g Ric$ is constant. This also means that $\delta Ric = 0$. In this case, we can conclude that $Ric \in \mathrm{H}(S^2 M)$. Moreover, if $(M, g)$ is a Hadamard manifold and $\|Ric\| \in L^q(M)$ at least for one $q \in (0, +\infty)$ and $q \neq 1$, then $Ric \equiv 0$ by Theorem 2.1. Next, we need to prove one obvious statement: If the sectional curvature is non-positive and the Ricci curvature is zero, then the Riemannian manifold is flat.

Namely, let $X \in T_x M$ be a unite vector and we complete it to an orthonormal basis $\{X, e_2, \ldots, e_n\}$ for $T_x M$ at an arbitrary point $x \in M$, then (see [1, p. 86])

$$Ric(X, X) = \sum_{a=2}^{n} \sec(X \wedge e_a).$$

In this case, from the conditions $Ric \equiv 0$ and $\sec \leq 0$ we obtain $\sec \equiv 0$, i.e., the sectional curvature vanishes identically. In this case, $(M, g)$ is a flat manifold. Again $(M, g)$ is simply connected, hence it follows that $(M, g)$ is isometric to the Euclidean space $\mathbf{R}^n$.

**Remark**. The author of the article [28] gives in it a complete classification of all 3-dimensional simply connected and complete $A$-spaces. This classification includes Riemannian globally symmetric spaces which are a well-known example of Hadamard manifolds.

## 4. Symmetric Killing tensors and Killing vectors on Hadamard manifolds

A tensor field $\varphi \in C^\infty(S^p M)$ satisfying the equation $\delta^* \varphi = 0$ is well-known in general relativity as a symmetric Killing $p$-tensor (see [30, p. 559]; [31]). It is the natural generalization of a *Killing vector* or, in other words, *an infinitesimal isometric* (see [6, pp. 42-43]; [30, p. 292];). Killing vector and Killing symmetric tensors have numerous applications. For example, while Killing vectors give the linear first integrals of the geodesic equations, Killing tensors give the quadratic, cubic, and higher-order first integrals on Riemannian manifolds.

We recall here the result of Dairbekov and Sharafutdinov in [32] and Heil in [33], respectively. Namely, they proved the following theorem: On a compact Riemannian manifold $(M, g)$ with non-positive sectional curvature, every Killing tensor is parallel. If there is a point in $(M, g)$ at which the sectional curvature is negative on all two-dimensional planes, then it is proportional to a power of the metric.

In this section, we consider symmetric Killing tensors on Hadamard manifolds and complement the above result.

The Ricci tensor $Ric$ of an $A$-space is a symmetric traceless Killing 2-tensor and, moreover, is a $\Delta_S$-harmonic tensor. We will consider a generalization of this concept to the case of a symmetric traceless Killing $p$-tensor for all $p \geq 2$.

Namely, let a smooth symmetric Killing $p$-tensor $\varphi$ is a traceless, then from the equation $\delta^*\varphi = 0$ we obtain $\delta\varphi = 0$. In this case, $\varphi \in \mathrm{H}(S^p M)$. Moreover, we proved the following theorem: Let $\varphi \in C^\infty(S^p M)$ be a divergence-free symmetric Killing tensor on a Riemannian manifold $(M, g)$, then it satisfies the following systems of differential equations $\Delta_S \varphi = 0$ and $\delta\varphi = 0$. Conversely, if $(M, g)$ is compact and $\varphi \in C^\infty(S^p M)$ satisfy the equations $\Delta_S \varphi = 0$ and $\delta\varphi = 0$, then $\varphi$ is a divergence-free symmetric Killing tensor (see [22]). In particular, this theorem is a natural generalization of the classical theorem on Killing vectors (see [6, p. 44]; [29, p. 44]).

Next, if $(M, g)$ is a Hadamard manifold and $\|\varphi\| \in L^q(M)$ at least for one $q \in (0, +\infty)$ and $q \neq 1$ for a divergence-free symmetric Killing tensor $\varphi$, then $\varphi \equiv 0$ by Theorem 2.1. In this case, we can formulate a theorem.

**Theorem 3.1.** *The Hadamard manifold $(M, g)$ does not admit a non-zero symmetric Killing $p$-tensor $\varphi \in C^\infty(S_0^p M)$ such that $\|\varphi\| \in L^q(M)$ at least for one $q \in (0, \infty)$ and $q \neq 1$.*

On the other hand, using Theorem 2.3, we can formulate a theorem analogous to Theorem 3.1.

**Theorem 3.2.** *The Hadamard manifold $(M, g)$ does not admit a non-zero symmetric Killing $p$-tensor $\varphi \in C^\infty(S_0^p M)$ such that the square of its norm belongs to $L^q(M)$ at least for one $q \in (0, \infty)$.*

The equation $\delta^*\varphi = 0$ for $\varphi = g(V, \cdot)$ defines the Killing vector $V$ (see [6, pp. 42 - 43]). From the equation $\delta^*\varphi = 0$ one can obtain (see also [6, p. 56])

$$\frac{1}{2}\left(\nabla^2 \|V\|^2\right)(X,X) = \|\nabla_X V\|^2 - g(R(V,X)X,V)$$

for an arbitrary $X \in TM$. If, in addition, $(M,g)$ is a Hadamard manifold, then $g(R(V,X)X,V) \leq 0$, and hence $f = \|V\|^2$ is a non-negative convex function. In this case, based of Corollary 1.4, we can formulate the following proposition.

**Corollary 3.3**. *The Hadamard manifold $(M,g)$ does not admit a non-zero Killing vector $V$ such that $f = \|V\|^2$ is a $L^q$-function at least for one $q \in (0, \infty)$.*

**Remark.** Corollary 3.3. generalizes the classical theorem on Killing vector fields on compact manifolds (see [29, p. 44]) and complements the following theorem: If $(M,g)$ is a complete Riemannian manifold with non-positive Ricci curvature, then every Killing vector on $(M,g)$ with finite global norm is parallel (see [35]).

## 5. Killing-Yano tensors on Riemannian globally symmetric spaces of noncompact type

The notion of *Killing-Yano tensors* were introduced into physics by Penrose et al. (see [36]; [37]; and etc.), and they played an important role in the development of general relativity (see, for example, [30, pp. 559-563]). On the other hand, these tensors were introduced in differential geometry thanks to K. Yano (see [29, p. 68]; [38]) and were fruitfully studied for a long time in the geometry of Riemannian manifolds (see, for example, [32]; [34]; [39]; [20] and etc.). In particular, in [39] was proved that a compact simply connected symmetric space carries a non-parallel Killing-Yano $p$-tensor $(p \geq 2)$ if and only if it isometric to a Riemannian product $\mathbf{S}^k \times N$, where $\mathbf{S}^k$ is a round sphere and $k > p$. In turn, we give here some applications of the results on subharmonic functions to the geometry of Killing-Yano $p$-tensors on symmetric spaces of non-compact type.

Let $(M, g)$ be an $n$-dimensional Riemannian manifold. It is obvious that the following equality holds

$$\dim \Lambda^p(T_x^*M) = \binom{n}{p}$$

for the vector space $\Lambda^p(T_x^*M)$ of covariant skew-symmetric $p$-tensors $(1 \leq p \leq n-1)$ on $T_xM$ at an arbitrary point $x \in M$.

A Killing-Yano $p$-tensor (or, in other words, Killing $p$-form) on $(M, g)$ is a skew-symmetric tensor whose covariant derivative is totally skew-symmetric, i.e., by definition, if $\omega \in C^\infty(\Lambda^p M)$ is a Killing-Yano tensor, then $\nabla \omega \in C^\infty(\Lambda^{p+1} M)$.

Let $d : C^\infty(\Lambda^p M) \to C^\infty(\Lambda^{p+1} M)$ be the operator of exterior derivative and $\delta : C^\infty(\Lambda^p M) \to C^\infty(\Lambda^{p-1} M)$ be the codifferentiation operator which defined as the canonical formal adjoint of $d$ (see [1, pp. 334-335]). Using these operators, one constructs the well-known Hodge – de Rham Laplacian $\Delta_H = \delta d + d \delta$ which admits a Weitzenböck decomposition (see [1, p. 347]; [17, pp. 77; 79]; [20])

$$\Delta_H \omega = \overline{\Delta} \omega + \Re(\omega) \qquad (5.1)$$

for any $\omega \in C^\infty(\Lambda^p M)$ and an algebraic symmetric operator $\Re : \Lambda^p M \to \Lambda^p M$ that is the Weitzenböck curvature operator of the ordinary Lichnerowicz Laplacian which was restricted to skew-symmetric $p$-tensors.

We define a non-negative scalar function by the equality $f = \|\omega\|$. Then using (5.1) we can write the well-known Bochner-Weizenbeck formula (see also [20])

$$\frac{1}{2} \Delta f^2 = -g(\Delta_H \omega, \omega) + \|\nabla \omega\|^2 + g(\Re(\omega), \omega). \qquad (5.2)$$

An arbitrary Killing-Yano tensor $\omega \in C^\infty(\Lambda^p M)$ satisfies the equation (see [20])

$$\Delta_H \omega = \frac{p+1}{p} \Re(\omega). \qquad (5.3)$$

In this case, from (5.2) and (5.3) we obtain the inequality

$$\frac{1}{2}\Delta f^2 \geq -\frac{1}{p} g(\Re(\omega),\omega). \tag{5.4}$$

We must recall that the Riemannian curvature tensor $R$ of $(M,g)$ defines a symmetric algebraic operator $\bar{R}: \Lambda^2(T_xM) \to \Lambda^2(T_xM)$ on the vector space $\Lambda^2(T_xM)$ of 2-forms over tangent space $T_xM$ at an arbitrary point $x \in M$ (see [1, pp. 82-83]; [17, p. 51]). There have been many papers on the relationship between the behavior of the curvature operator $\bar{R}$ of a Riemannian manifold $(M,g)$ and some global characterization of it, such as its homotopy type, topological types and etc. We say that the manifold $(M,g)$ has a *non-positive curvature operator* $\bar{R}$ if $g(\bar{R}(\theta),\theta) \leq 0$ for all nonzero two-forms $\theta \neq 0$. At the same time, it can be also concluded here that if the curvature operator $\bar{R}$ is non-positive then the quadratic form $g(\Re(\omega),\omega) \leq 0$ for any $\omega \in \Lambda^p M$ by the formulas in [1 pp. 345-346] and [2, p. 29-30]. Side by side, we know that a Riemannian symmetric space has non-positive curvature operator $\bar{R}$ if and only if it has non-positive sectional curvature (see [24]). In this case, we can conclude from (5.4) that

$$\Delta f^2 \geq 0 \tag{5.5}$$

for $f = \|\omega\|$. Hence $f^2 = \|\omega\|^2$ is a subharmonic function. At the same time, we known that a Riemannian globally symmetric space of noncompact type does not admit a nonzero non-negative subharmonic $L^q$-function for $q \in (0,\infty)$. Therefore, we can formulate the following proposition.

**Theorem 5.1**. *An $n$-dimensional Riemannian globally symmetric space of noncompact type does not admit a non-zero Killing-Yano $p$-tensor $(1 \leq p \leq n-1)$ such that the square of its norm is an $L^q$-function at least one $q \in (0,+\infty)$.*

The following equality holds

$$\frac{1}{2}\Delta f^2 = f \cdot \Delta f + \|\nabla f\|^2 \tag{5.6}$$

for $f = \|\omega\|$. By virtue of the Kato inequality (see [7, p. 380])

$$\|\nabla\|\omega\|\| \leq \|\nabla \omega\|$$

and (5.6) the Bochner-Weitzenböck formula (5.2) can be rewritten as the inequality

$$f \Delta f \geq -g(\Delta_H \omega, \omega) + g(\Re(\omega), \omega). \tag{5.7}$$

In turn, a Killing-Yano tensor $\omega \in C^\infty(\Lambda^p M)$ satisfies the equation (see [20])

$$\Delta_H \omega = \frac{p+1}{p}\Re(\omega). \tag{5.8}$$

Then from (5.7) and (5.8) we obtain the inequality

$$f \Delta f \geq -\frac{1}{p} g(\Re(\omega), \omega). \tag{5.9}$$

Next we will use Lemma 1.3. Namely, if $f = \|\omega\|$ is a $L^q$-function for some $q \neq 1$, then $\|\omega\| = constant$. This means that the identity $\omega \equiv 0$ holds on a Hadamard manifold and, in particular, on a Riemannian globally symmetric space of noncompact type. Therefore, we can formulate the following proposition.

**Theorem 5.2**. *An $n-$ dimensional Riemannian globally symmetric space of noncompact type does not admit a non-zero Killing-Yano $p-$tensor $(1 \leq p \leq n-1)$ such that its norm is a $L^q$-function at least for one $q \in (0, +\infty)$ and $q \neq 1$.*

For any $p$ with $0 \leq p \leq n$, we can define the Hodge operator $*$ to be the unique vector-bundle isomorphism (see [17, p. 33])

$$*: \Lambda^p M \to \Lambda^{n-p} M$$

such that $*^2 = (-1)^{p(n-p)}$ and

$$\omega \wedge (*\omega') = g(\omega, \omega') dvol_g$$

for any $\omega, \omega' \in \Lambda^p M$ and the volume form $dvol_g$ of $(M, g)$. Side by side, $p$-tensor $*\omega$ for an arbitrary Killing-Yano $(n-p)$-tensor $\omega$ is called *closed conformal Killing-Yano $p$-tensor* or, in other words, *closed conformal Killing $p$-form* (see [41]; [42]). In particular, for any closed conformal Killing-Yano tensor $\omega \in C^\infty(\Lambda^p M)$ we have (see [41])

$$\Delta_H \omega = \frac{n-p+1}{n-p} \Re(\omega). \tag{5.10}$$

Further, a direct calculation based on (5.7) and (5.10) we obtain the inequality

$$\frac{1}{2}\Delta f^2 \geq -\frac{1}{n-p} g(\Re(\omega), \omega) \tag{5.11}$$

for a closed conformal Killing-Yano $p$-tensor $\omega$. Then arguments similar to those carried out above for the case of the Killing-Yano $p$-tensor allow us to formulate the following theorem.

**Theorem 5.3**. *An $n$-dimensional Riemannian globally symmetric space of noncompact type does not admit a non-zero closed conformal Killing-Yano $p$-tensor $(1 \leq p \leq n-1)$ such that the square of its norm is an $L^q$-function at least one $q \in (0, +\infty)$.*

Then arguments similar to those carried out above for the case of the Killing-Yano $p$-tensor allow us to formulate the following theorem.

**Theorem 5.4**. *An $n$-dimensional Riemannian globally symmetric space of noncompact type does not admit a non-zero closed conformal Killing-Yano $p$-*

*tensor $(1 \leq p \leq n-1)$ such that its norm is a $L^q$-function at least for one $q \in (0, +\infty)$ and $q \neq 1$.*

**Remark.** The theorems of this section generalize the classical theorems on Killing-Yano and closed conformal Killing-Yano tensors on compact manifolds (see [29, p. 68-70]; [41]) and complement our results in [42].